\documentclass[12pt]{amsart}
\usepackage{t1enc}
\usepackage[latin1]{inputenc}
\usepackage{amsmath,amsfonts,amssymb,amsthm}
\usepackage{ulem,xcolor,cancel}

\newtheorem{theorem}{Theorem}
\newtheorem{lemma}{Lemma}

\newtheorem{corollary}{Corollary}

\newtheorem{remark}{Remark}
\numberwithin{equation}{subsection}

\begin{document}
\author{I. Blahota, L. E. Persson, G. Tephnadze}
\title[N\"orlund means]{On the N\"orlund means of Vilenkin-Fourier series }
\address{I. Blahota, Institute of Mathematics and Computer Sciences, College
of Ny\'\i regyh\'aza, P.O. Box 166, Ny\'\i regyh\'aza, H-4400, Hungary.}
\email{blahota@nyf.hu}
\address{L.-E. Persson, Department of Engineering Sciences and Mathematics,
Luleå\ University of Technology, SE-971 87 Luleå, Sweden and Narvik
University College, P.O. Box 385, N-8505, Narvik, Norway.}
\email{larserik@ltu.se}
\address{G. Tephnadze, Department of Mathematics, Faculty of Exact and
Natural Sciences, Tbilisi State University, Chavchavadze str. 1, Tbilisi
0128, Georgia and Department of Engineering Sciences and Mathematics, Luleå\
University of Technology, SE-971 87 Luleå, Sweden.}
\email{giorgitephnadze@gmail.com}
\thanks{Supported by T\'AMOP 4.2.2.A-11/1/KONV-2012-0051 and by Shota
Rustaveli National Science Foundation grant no. 52/54 (Bounded operators on
the martingale Hardy spaces).}
\date{}

\begin{abstract}
In this paper we prove and discuss some new $\left( H_{p},L_{p}\right)$-type
inequalities of weighted maximal operators of Vilenkin-N\"orlund means with
non-increasing coefficients. These results are the best possible in a
special sense. As applications, both some well-known and new results are
pointed out in the theory of strong convergence of Vilenkin-N\"orlund means
with non-increasing coefficients.
\end{abstract}

\maketitle



\section{Introduction}

\bigskip The definitions and notations used in this introduction can be
found in our next Section. In the one-dimensional case the weak (1,1)-type
inequality for maximal operator of Fej\'er means $\sigma ^{\ast }$ can be
found in Schipp \cite{Sc} for Walsh series and in P\'al, Simon \cite{PS} for
bounded Vilenkin series. Fujji \cite{Fu} and Simon \cite{Si2} verified that $%
\sigma ^{\ast }$ is bounded from $H_{1}$ to $L_{1}$. Weisz \cite{We2}
generalized this result and proved boundedness of $\sigma ^{\ast }$ from the
martingale space $H_{p}$ to the Lebesgue space $L_{p}$ for $p>1/2$. Simon 
\cite{Si1} gave a counterexample, which shows that boundedness does not hold
for $0<p<1/2.$ A counterexample for $p=1/2$ was given by Goginava \cite{Go}.
Weisz \cite{We4} proved that the maximal operator of the Fej\'er means $%
\sigma ^{\ast }$ is bounded from the Hardy space $H_{1/2}$ to the space $%
weak-L_{1/2}$. Goginava \cite{GoSzeged} (see also \cite{tep2}) proved that
weighted maximal operator $\widetilde{\sigma }^{\ast }$ is bounded from the
Hardy space $H_{1/2}$ to the space $L_{1/2}$. Moreover, the rate of the
weights $\left\{ \log ^{2}\left( n+1\right) \right\} _{n=1}^{\infty }$ in $n$%
-th Fej\'er mean is given exactly. Analogical results for $0<p<1/2$ were
proved in \cite{tep3}.

Riesz's logarithmic means with respect to Walsh and Vilenkin systems were
studied by several authors. We mention, for instance, the papers by Simon 
\cite{Si1}, G\'at, Nagy \cite{GN}. In \cite{tep5} it was proved that the
maximal operator of Riesz's means $R^{\ast }$ is bounded from the Hardy
space $H_{1/2}$ to the space $weak-L_{1/2}$, but is not bounded from the
Hardy space $H_{p}$ to the space $L_{p},$ when $0<p\leq 1/2. $ Moreover,
there were proved some theorems of boundedness of weighted maximal operators
of Riesz's logarithmic means, with respect to Vilenkin-Fourier series.

Móricz and Siddiqi \cite{Mor} investigated the approximation properties of
some special N\"orlund means of Walsh-Fourier series of $L_{p}$ function in
norm. The case when $q_{k}=1/k$ was excluded, since the methods of Móricz
and Siddiqi are not applicable to N\"orlund logarithmic means. In \cite{Ga2}
G\'at and Goginava investigated some properties of the N\"orlund logarithmic
means of functions in the class of continuous functions and in the Lebesgue
space $L_{1}.$ In \cite{tep4} it was proved that there exists a martingale $%
f\in H_{p},\ (0<p\leq 1),$ such that the maximal operator of N\"orlund
logarithmic means $L^{\ast }$ is not bounded in the space $L_{p}.$ For more
information on N\"orlund logarithmic means, see paper of Blahota and G\'at 
\cite{bg} and Nagy (see \cite{n}, \cite{nagy} and \cite{nag}).

In \cite{gog8} Goginava investigated the behaviour of Cesàro means of
Walsh-Fourier series in detail. In the two-dimensional case approximation
properties of N\"orlund and Cesàro means was considered by Nagy \cite{nn}.
Weisz \cite{we6} proved that the maximal operator $\sigma ^{\alpha ,\ast }$
is bounded from the martingale space $H_{p}$ to the space $L_{p}$ for $%
p>1/\left( 1+\alpha \right) .$ Goginava \cite{gog4} gave a counterexample,
which shows that boundedness does not hold for $0<p\leq 1/\left( 1+\alpha
\right) .$ Simon and Weisz \cite{sw} showed that the maximal operator $%
\sigma ^{\alpha ,\ast }\ \left( 0<\alpha <1\right) $ of the $\left( C,\alpha
\right) $ means is bounded from the Hardy space $H_{1/\left( 1+\alpha
\right) }$ to the space $weak-L_{1/\left( 1+\alpha \right) }$. In \cite{bt}
it was also proved that the maximal operator $\widetilde{\sigma }^{\alpha
,\ast }$ is bounded from the Hardy space $H_{1/\left( 1+\alpha \right) }$ to
the space $L_{1/\left( 1+\alpha \right) }$. Moreover, this result can not be
improved in the following sense:

\textbf{Theorem BT.} \textit{(Blahota, Tephnadze \cite{bt}) Let $\mathit{%
0<\alpha \leq 1}$ and $\varphi :\mathbb{N}_{+}\rightarrow \lbrack 1,\infty )$
be a non-decreasing function satisfying the condition 
\begin{equation*}
\overline{\lim_{n\rightarrow \infty }}\frac{\log ^{1+\alpha }n}{\varphi
\left( n\right) }=\infty .
\end{equation*}%
Then there exists a martingale $f\in H_{1/(1+\alpha )}\left( G\right) ,$
such that 
\begin{equation*}
\sup_{n\in \mathbb{N}}\left\Vert \frac{\sigma _{n}^{\alpha }f}{\varphi
\left( n\right) }\right\Vert _{1/\left( 1+\alpha \right) }=\infty .
\end{equation*}%
}

It is well-known that Vilenkin systems do not form bases in the space $%
L_{1}\left( G_{m}\right) $. Moreover, there is a function in the Hardy space 
$H_{1}\left( G_{m}\right) $, such that the partial sums of $f$ are not
bounded in $L_{1}$-norm. However, in Gát \cite{Ga1} (see also \cite{b}) the
following strong convergence result was obtained for all $f\in H_{1}:$ 
\begin{equation*}
\underset{n\rightarrow \infty }{\lim }\frac{1}{\log n}\overset{n}{\underset{%
k=1}{\sum }}\frac{\left\Vert S_{k}f-f\right\Vert _{1}}{k}=0.
\end{equation*}

Simon \cite{Si4} (see also \cite{tep8}) proved that there exists an absolute
constant $c_{p},$ depending only on $p,$ such that 
\begin{equation}
\frac{1}{\log ^{\left[ p\right] }n}\overset{n}{\underset{k=1}{\sum }}\frac{%
\left\Vert S_{k}f\right\Vert _{p}^{p}}{k^{2-p}}\leq c_{p}\left\Vert
f\right\Vert _{H_{p}}^{p},\ \ \ \left( 0<p\leq 1\right)  \label{1cc}
\end{equation}%
for all $f\in H_{p}$ and $n\in \mathbb{N}_{+},$ where $\left[ p\right] $
denotes the integer part of $p.$ In \cite{tep9} it was proved that sequence $%
\left\{ 1/k^{2-p}\right\} _{k=1}^{\infty }\ \left( 0<p<1\right) $ in (\ref%
{1cc}) can not be improved.

Weisz considered the norm convergence of Fej\'er means of Vilenkin-Fourier
series and proved the following:

\textbf{Theorem W1.} \textit{(Weisz \cite{We1}) Let $p>1/2$ and $f\in H_{p}.$
Then there exists an absolute constant $c_{p}$, depending only on $p$, such
that 
\begin{equation*}
\left\Vert \sigma _{k}f\right\Vert _{p}\leq c_{p}\left\Vert f\right\Vert
_{H_{p}},\text{ for all }f\in H_{p}\text{ and }k=1,2,\dots.
\end{equation*}%
} Theorem W1 implies that 
\begin{equation*}
\frac{1}{n^{2p-1}}\overset{n}{\underset{k=1}{\sum }}\frac{\left\Vert \sigma
_{k}f\right\Vert _{p}^{p}}{k^{2-2p}}\leq c_{p}\left\Vert f\right\Vert
_{H_{p}}^{p},\ \ \ \left( 1/2<p<\infty,\ n=1,2,\dots \right).
\end{equation*}

If Theorem W1 holds for $0<p\leq 1/2,$ then we would have 
\begin{equation}
\frac{1}{\log ^{\left[ 1/2+p\right] }n}\overset{n}{\underset{k=1}{\sum }} 
\frac{\left\Vert \sigma _{k}f\right\Vert _{p}^{p}}{k^{2-2p}}\leq
c_{p}\left\Vert f\right\Vert _{H_{p}}^{p}, \ \ \ \left( 0<p\leq 1/2,\
n=2,3,\dots \right).  \label{2cc}
\end{equation}

However, in \cite{tep1} it was proved that the assumption $p>1/2$ in Theorem
W1 is essential. In particular, we showed that there exists a martingale $%
f\in H_{1/2},$ such that \ 
\begin{equation*}
\sup_{n}\left\Vert \sigma _{n}f\right\Vert _{1/2}=\infty .
\end{equation*}

In \cite{bt1} it was proved that (\ref{2cc}) holds, though Fej\'er means is
not of type $\left( H_{p},L_{p}\right),$ for $0<p\leq 1/2.$ This result for $%
\left( C,\alpha \right) \ \left( 0<\alpha <1\right) $ means when $p=1/\left(
1+\alpha \right) $ was generalized in \cite{bt}.

In this paper we prove and discuss some new $\left( H_{p},L_{p}\right)$-type
inequalities of weighted maximal operators of Vilenkin-N\"orlund means with
non-increasing coefficients. These results are the best possible in a
special sense. As applications, both some well-known and new results are
pointed out in the theory of strong convergence of Vilenkin-N\"orlund means.

This paper is organized as follows: in order not to disturb our discussions
later on some definitions and notations are presented in Section 2. The main
results and some of its consequences can be found in Section 3. For the
proofs of the main results we need some auxiliary results of independent
interest. Also these results are presented in Section 3. The detailed proofs
are given in Section 4.

\section{Definitions and Notations}

Denote by $\mathbb{N}_{+}$ the set of the positive integers, $\mathbb{N}:=%
\mathbb{N}_{+}\cup \{0\}.$ Let $m:=(m_{0},m_{1},\dots)$ be a sequence of the
positive integers not less than 2. Denote by 
\begin{equation*}
Z_{m_{n}}:=\{0,1,\ldots,m_{n}-1\}
\end{equation*}%
the additive group of integers modulo $m_{n}$.

Define the group $G_{m}$ as the complete direct product of the groups $%
Z_{m_{n}}$ with the product of the discrete topologies of $Z_{m_{n}}`$s. In
this paper we discuss bounded Vilenkin groups, i.e. the case when $\sup_{n\in%
\mathbb{N}}m_{n}<\infty.$

The direct product $\mu $ of the measures 
\begin{equation*}
\mu _{n}\left( \{j\}\right) :=1/m_{n},\ (j\in Z_{m_{n}})
\end{equation*}%
is the Haar measure on $G_{m}$ with $\mu \left( G_{m}\right) =1.$

The elements of $G_{m}$ are represented by sequences 
\begin{equation*}
x:=\left( x_{0},x_{1},\ldots,x_{n},\ldots \right),\ \left( x_{n}\in
Z_{m_{n}}\right).
\end{equation*}

It is easy to give a base for the neighbourhood of $G_{m}:$

\begin{equation*}
I_{0}\left( x\right) :=G_{m},\ I_{n}(x):=\{y\in G_{m}\mid
y_{0}=x_{0},\ldots,y_{n-1}=x_{n-1}\}\,\,\left( x\in G_{m},\text{ }n\in 
\mathbb{N}\right).
\end{equation*}

Denote $I_{n}:=I_{n}\left( 0\right),$ for $n\in \mathbb{N}_{+}$ and 
\begin{equation*}
e_{n}:=\left( 0,\ldots,x_{n}=1,0,\ldots \right) \in G_{m},\ \left( n\in 
\mathbb{N}\right).
\end{equation*}

It is evident that 
\begin{equation}
\overline{I_{N}}=\left( \overset{N-2}{\underset{k=0}{\bigcup }}\overset{
m_{k}-1}{\underset{x_{k}=1}{\bigcup }}\overset{N-1}{\underset{l=k+1}{\bigcup 
}}\overset{m_{l}-1}{\underset{x_{l}=1}{\bigcup }}I_{l+1}\left(
x_{k}e_{k}+x_{l}e_{l}\right) \right) \bigcup \left(\underset{k=0}{
\bigcup\limits^{N-1}}\overset{m_{k}-1}{\underset{x_{k}=1}{\bigcup }}
I_{N}\left( x_{k}e_{k}\right) \right).  \label{2}
\end{equation}

\bigskip If we define the so-called generalized number system based on $m$
in the following way : 
\begin{equation*}
M_{0}:=1,\ M_{n+1}:=m_{n}M_{n}\ \ \ (n\in \mathbb{N}),
\end{equation*}%
then every $n\in \mathbb{N}$ can be uniquely expressed as 
\begin{equation*}
n=\sum_{k=0}^{\infty }n_{k}M_{k},\text{ where }n_{k}\in Z_{m_{k}} \ (k\in 
\mathbb{N}_{+})
\end{equation*}
and only a finite number of $n_{k}`$s differ from zero.

Next, we introduce on $G_{m}$ an orthonormal system which is called the
Vilenkin system. At first we define the complex-valued function $r_{k}\left(
x\right) :G_{m}\rightarrow\mathbb{C},$ the generalized Rademacher functions,
by 
\begin{equation*}
r_{k}\left( x\right) :=\exp \left( 2\pi ix_{k}/m_{k}\right),\text{ }\left(
i^{2}=-1,x\in G_{m},\text{ }k\in\mathbb{N}\right).
\end{equation*}

Now, define the Vilenkin system $\,\,\,\psi :=(\psi _{n}:n\in\mathbb{N})$ on 
$G_{m}$ as: 
\begin{equation*}
\psi _{n}(x):=\prod\limits_{k=0}^{\infty }r_{k}^{n_{k}}\left( x\right) \ \ \
\left( n\in\mathbb{N}\right).
\end{equation*}

Specifically, we call this system the Walsh-Paley system, when $m\equiv 2.$

The norm (or quasi-norm) of the space $L_{p}(G_{m})\ \left( 0<p<\infty
\right) $ is defined by 
\begin{equation*}
\left\Vert f\right\Vert _{p}^{p}:=\int_{G_{m}}\left\vert f\right\vert
^{p}d\mu.
\end{equation*}
\qquad

The space $weak-L_{p}\left( G_{m}\right) $ consists of all measurable
functions $f$, for which 
\begin{equation*}
\left\Vert f\right\Vert _{weak-L_{p}}^{p}:=\underset{\lambda >0}{\sup}
\,\lambda^{p}\mu \left( f>\lambda \right) <\infty.
\end{equation*}

The Vilenkin system is orthonormal and complete in $L_{2}\left( G_{m}\right) 
$ (see \cite{Vi}).

Now we introduce analogues of the usual definitions in Fourier-analysis. If $%
f\in L_{1}\left( G_{m}\right) $ we can define the Fourier coefficients, the
partial sums of the Fourier series, the Dirichlet kernels with respect to
the Vilenkin system in the usual manner:

\begin{equation*}
\widehat{f}\left( n\right) :=\int_{G_{m}}f\overline{\psi }_{n}d\mu, \ \
S_{n}f:=\sum_{k=0}^{n-1}\widehat{f}\left( k\right) \psi _{k}, \ \
D_{n}:=\sum_{k=0}^{n-1}\psi _{k},\ \ \left(n\in\mathbb{N}_{+}\right)
\end{equation*}
respectively.

Recall that 
\begin{equation}
D_{M_{n}}\left( x\right) =\left\{ 
\begin{array}{ll}
M_{n}, & \text{if\thinspace \thinspace \thinspace }x\in I_{n}, \\ 
0, & \text{if}\,\,x\notin I_{n}.%
\end{array}%
\right.  \label{1dn}
\end{equation}%
It is also known that (see \cite{blahota}, \cite{gg1} and \cite{go1}) 
\begin{equation}
D_{sM_{n}}=D_{M_{n}}\sum_{k=0}^{s-1}\psi
_{kM_{n}}=D_{M_{n}}\sum_{k=0}^{s-1}r_{n}^{k},  \label{2dn}
\end{equation}%
and 
\begin{equation}
D_{sM_{n}-j}=D_{sM_{n}}-w_{sM_{n}-1}\overline{D_{j}},\ j=1,\dots,M_{n}-1.
\label{3dn}
\end{equation}

The $\sigma $-algebra generated by the intervals $\left\{ I_{n}\left(
x\right) :x\in G_{m}\right\} $ will be denoted by $\digamma _{n}\left(n\in%
\mathbb{N} \right).$ Denote by $f=\left( f^{\left( n\right) },n\in \mathbb{N}%
\right)$ a martingale with respect to $\digamma _{n}\left( n\in \mathbb{N}%
\right).$ (for details see e.g. \cite{We1}).

The maximal function of a martingale $f$ is defined by 
\begin{equation*}
f^{\ast }:=\sup_{n\in\mathbb{N}}\left\vert f^{(n)}\right\vert.
\end{equation*}

For $0<p<\infty $ the Hardy martingale spaces $H_{p}\ \left( G_{m}\right) $
consist of all martingales, for which 
\begin{equation*}
\left\Vert f\right\Vert _{H_{p}}:=\left\Vert f^{\ast }\right\Vert
_{p}<\infty.
\end{equation*}

If $f=\left( f^{\left( n\right) },n\in\mathbb{N}\right) $ is a martingale,
then the Vilenkin-Fourier coefficients must be defined in a slightly
different manner: 
\begin{equation*}
\widehat{f}\left( i\right) :=\lim_{k\rightarrow \infty
}\int_{G_{m}}f^{\left( k\right) }\overline{\psi }_{i}d\mu.
\end{equation*}

Let $\{q_{n}:n\geq 0\}$ be a sequence of non-negative numbers. The $n$-th
N\"orlund mean is defined by

\begin{equation*}
t_{n}{f}:=\frac{1}{Q_{n}}\overset{n}{\underset{k=1}{\sum }}q_{n-k}S_{k}f,
\end{equation*}%
where \ 
\begin{equation*}
Q_{n}:=\sum_{k=0}^{n-1}q_{k}.
\end{equation*}

It is well known that 
\begin{equation*}
t_{n}f\left( x\right) =\int_{G_{m}}f\left( t\right) F_{n}\left( x-t\right)
dt,
\end{equation*}%
where $F_{n}$ are the so called N\"orlund kernels 
\begin{equation*}
F_{n}:=\frac{1}{Q_{n}}\overset{n}{\underset{k=1}{\sum }}q_{n-k}D_{k}.
\end{equation*}

We always assume that $q_{0}>0$ and $\lim_{n\rightarrow \infty }Q_{n}=\infty$%
. In this case (see \cite{moo}) the summability method generated by $%
\{q_{n}:n\geq 0\}$ is regular if and only if 
\begin{equation*}
\lim_{n\rightarrow \infty }\frac{q_{n-1}}{Q_{n}}=\infty.
\end{equation*}

If $q_{n}\equiv 1,$ then we get the usual $n$-th Fej\'er mean and Fej\'er
kernel 
\begin{equation*}
\sigma _{n}f:=\frac{1}{n}\sum_{k=1}^{n}S_{k}f, \ \ K_{n}:=\frac{1}{n}%
\sum_{k=1}^{n}D_{k},
\end{equation*}
respectively.

Let $t,n\in\mathbb{N}$. It is known that (see \cite{gat}) 
\begin{equation}
K_{M_{n}}\left( x\right) =\left\{ 
\begin{array}{ll}
0, & \text{ if }x-x_{t}e_{t}\notin I_{n},\ x\in I_{t}\backslash I_{t+1}, \\ 
\frac{M_{t}}{1-r_{t}\left( x\right) }, & \text{ if }x-x_{t}e_{t}\in I_{n}, \
x\in I_{t}\backslash I_{t+1}, \\ 
\left( M_{n}+1\right) /2, & \text{ if }x\in I_{n}.%
\end{array}
\right.  \label{5a}
\end{equation}

The $\left( C,\alpha \right) $-means (Cesàro means) of the Vilenkin-Fourier
series are defined by 
\begin{equation*}
\sigma _{n}^{\alpha }f=\frac{1}{A_{n}^{\alpha }}\overset{n}{\underset{k=1}{%
\sum }}A_{n-k}^{\alpha -1}S_{k}f,
\end{equation*}%
where \ 
\begin{equation*}
A_{0}^{\alpha }:=0,\ \ A_{n}^{\alpha }:=\frac{\left( \alpha +1\right) \ldots
\left( \alpha +n\right) }{n!},\ \alpha \neq -1,-2,\dots
\end{equation*}

The $n$-th Riesz logarithmic mean $R_{n}$ and the N\"orlund logarithmic mean 
$L_{n}$ are defined by 
\begin{equation*}
R_{n}f:=\frac{1}{l_{n}}\sum_{k=1}^{n-1}\frac{S_{k}f}{k},\ \ L_{n}f:=\frac{1}{%
l_{n}}\sum_{k=1}^{n-1}\frac{S_{k}f}{n-k},
\end{equation*}
respectively, where 
\begin{equation*}
l_{n}:=\sum_{k=1}^{n-1}1/k.
\end{equation*}

For the martingale $f$ we consider the following maximal operators: 
\begin{equation*}
t^{\ast }f:=\sup_{n\in \mathbb{N}}\left\vert t_{n}f\right\vert ,\ \ 
\end{equation*}%
\begin{equation*}
\sigma ^{\ast }f:=\sup_{n\in \mathbb{N}}\left\vert \sigma _{n}f\right\vert
,\ \ \ \sigma ^{\alpha ,\ast }f:=\sup_{n\in \mathbb{N}}\left\vert \sigma
_{n}^{\alpha }f\right\vert ,
\end{equation*}%
\begin{equation*}
R^{\ast }f:=\sup_{n\in \mathbb{N}}\left\vert R_{n}f\right\vert ,\ \ \ \
L^{\ast }f:=\sup_{n\in \mathbb{N}}\left\vert L_{n}f\right\vert .
\end{equation*}

We also define the following weighted maximal operators:%
\begin{equation*}
\widetilde{t}^{\ast }f:=\sup_{n\in \mathbb{N}}\left\vert t_{n}f\right\vert
/\log ^{1+\alpha }\left( n+1\right) ,\ \ \ 
\end{equation*}%
\begin{equation*}
\widetilde{\sigma }^{\alpha ,\ast }f:=\sup_{n\in \mathbb{N}}\left\vert
\sigma _{n}^{\alpha }f\right\vert /\log ^{1+\alpha }\left( n+1\right) ,
\end{equation*}%
\begin{equation*}
\widetilde{\sigma }^{\ast }f:=\sup_{n\in \mathbb{N}}\left\vert \sigma
_{n}f\right\vert /\log ^{2}\left( n+1\right) .
\end{equation*}

A bounded measurable function $a$ is called a $p$-atom, if there exists an
interval $I$, such that 
\begin{equation*}
\int_{I}ad\mu =0,\ \ \left\Vert a\right\Vert _{\infty }\leq \mu \left(
I\right) ^{-1/p},\text{ \ \ supp}\left( a\right) \subset I.
\end{equation*}

\section{Results}

\vspace{0.5cm}

\begin{center}
\textbf{Main results and some of its consequences}
\end{center}

\vspace{0.5cm}

\begin{theorem}
\label{Theorem1}Let $f\in H_{1/(1+\alpha )},$ where $0<\alpha \leq 1$ and $%
\{q_{n}:n\geq 0\},$ be a sequence of non-increasing numbers, such that 
\begin{equation}
n^{\alpha }/Q_{n}=O\left( 1\right) ,\text{ as }\ n\rightarrow \infty ,
\label{6a}
\end{equation}%
and%
\begin{equation}
\left( q_{n}-q_{n+1}\right) /n^{\alpha -2}=O\left( 1\right) ,\text{ as }\
n\rightarrow \infty .  \label{7a}
\end{equation}%
Then there exists an absolute constant $c_{\alpha },$ depending only on $%
\alpha ,$ such that 
\begin{equation*}
\left\Vert \overset{\sim }{t}^{\ast }f\right\Vert _{1/\left( 1+\alpha
\right) }\leq c_{\alpha }\left\Vert f\right\Vert _{H_{1/\left( 1+\alpha
\right) }}.
\end{equation*}
\end{theorem}

\begin{corollary}
\label{Corollary1} (Blahota, Tephnadze \cite{bt}) Let $f\in H_{1/(1+\alpha
)},$ where $0<\alpha <1.$ Then there exists an absolute constant $c_{\alpha
},$ depending only on $\alpha,$ such that 
\begin{equation*}
\left\Vert \overset{\sim }{\sigma }^{\alpha,\ast }f\right\Vert _{1/\left(
1+\alpha \right) }\leq c_{\alpha }\left\Vert f\right\Vert _{H_{1/\left(
1+\alpha \right) }}.
\end{equation*}
\end{corollary}

\begin{corollary}
\label{Corollary2} (Goginava \cite{GoSzeged}, Tephnadze \cite{tep2}) Let $%
f\in H_{1/2}.$ Then there exists an absolute constant $c$, such that 
\begin{equation*}
\left\Vert \overset{\sim }{\sigma }^{\ast }f\right\Vert _{1/2}\leq
c\left\Vert f\right\Vert _{H_{1/2}}.
\end{equation*}
\end{corollary}

\begin{theorem}
\label{Theorem2} Let $f\in H_{1/(1+\alpha )},$ where $0<\alpha <1$ and $%
\{q_{n}:n\geq 0\},$ be a sequence of non-increasing numbers, satisfying
condition (\ref{6a}) and (\ref{7a}). Then there exists an absolute constant $%
c_{\alpha },$ depending only on $\alpha,$ such that 
\begin{equation*}
\frac{1}{\log n}\overset{n}{\underset{k=1}{\sum }}\frac{\left\Vert
t_{k}f\right\Vert _{H_{1/\left( 1+\alpha \right) }}^{1/\left( 1+\alpha
\right) }}{k}\leq c_{\alpha }\left\Vert f\right\Vert _{H_{1/\left( 1+\alpha
\right) }}^{1/\left( 1+\alpha \right) }.
\end{equation*}
\end{theorem}

\begin{corollary}
\label{Corollary3}(Blahota, Tephnadze \cite{bt}) Let $f\in H_{1/(1+\alpha
)}, $ where $0<\alpha <1.$ Then there exists an absolute constant $c_{\alpha
},$ depending only on $\alpha,$ such that 
\begin{equation*}
\frac{1}{\log n}\overset{n}{\underset{k=1}{\sum }}\frac{\left\Vert \sigma
_{k}^{\alpha }f\right\Vert _{H_{1/\left( 1+\alpha \right) }}^{1/\left(
1+\alpha \right) }}{k}\leq c_{\alpha }\left\Vert f\right\Vert _{H_{1/\left(
1+\alpha \right) }}^{1/\left( 1+\alpha \right) }.
\end{equation*}
\end{corollary}

\begin{corollary}
\label{Corollary4}(Blahota, Tephnadze \cite{bt1}, Tephnadze \cite{tep6}) Let 
$f\in H_{1/2}.$ Then there exists an absolute constant $c,$ such that 
\begin{equation*}
\frac{1}{\log n}\overset{n}{\underset{k=1}{\sum }}\frac{\left\Vert \sigma
_{k}f\right\Vert _{1/2}^{1/2}}{k}\leq c\left\Vert f\right\Vert
_{H_{1/2}}^{1/2}.
\end{equation*}
\end{corollary}

\begin{remark}
\label{Remark1} For some $\{q_{n}:n\geq 0\}$ sequences of non-increasing
numbers conditions (\ref{6a}) and (\ref{7a}) can be true or false
independently.
\end{remark}

\begin{remark}
\label{Remark2}Since Cesàro means satisfy conditions (\ref{6a}) and (\ref{7a}%
), we immediately obtain from the Theorem BT that the rate of the weights $%
\left\{ \log ^{1+\alpha }\left( n+1\right) \right\} _{n=1}^{\infty }$ in $n$%
-th Nörlund mean can not be improved.
\end{remark}

\vspace{1cm}

\begin{center}
\textbf{Some auxiliary results}
\end{center}

\vspace{0.5cm}

Weisz proved that the following is true:

\begin{lemma}
(Weisz \cite{We1}) \label{W}Suppose that an operator $T$ is $\sigma $-linear
and for some $0<p\leq 1$ 
\begin{equation*}
\int\limits_{\overline{I}}\left\vert Ta\right\vert ^{p}d\mu \leq
c_{p}<\infty,
\end{equation*}%
for every $p$-atom $a$, where $I$ denotes the support of the atom. If $T$ is
bounded from $L_{\infty \text{ }}$ to $L_{\infty },$ then 
\begin{equation*}
\left\Vert Tf\right\Vert _{p}\leq c_{p}\left\Vert f\right\Vert _{H_{p}}.
\end{equation*}
\end{lemma}

We also state three new Lemmas we need for the proofs of our main results
but which are also of independent interest:

\begin{lemma}
\label{Lemma2}Let $sM_{n}<r\leq \left( s+1\right) M_{n},$ where $1\leq
s<m_{n}.$ Then 
\begin{equation}
Q_{r}F_{r}=Q_{r}D_{sM_{n}}-w_{sM_{n}-1}\overset{sM_{n}-2}{\underset{l=1}{%
\sum }}\left( q_{r-sM_{n}+l}-q_{r-sM_{n}+l+1}\right) l\overline{K_{l}}
\label{8a}
\end{equation}%
\begin{equation*}
-w_{sM_{n}-1}\left( sM_{n}-1\right) q_{r-1}\overline{K_{sM_{n}-1}}%
+w_{sM_{n}}Q_{r-sM_{n}}F_{r-sM_{n}}.
\end{equation*}
\end{lemma}

The next Lemma is generalization of analogical estimation of Cesàro means
(see \cite{gago})

\begin{lemma}
\label{Lemma3}Let $0<\alpha \leq 1$ and $\{q_{n}:n\geq 0\}$ be a sequence of
non-increasing numbers, satisfying conditions (\ref{6a}) and (\ref{7a}). 
\textit{Then} 
\begin{equation*}
\left\vert F_{n}\right\vert \leq \frac{c_{\alpha }}{n^{\alpha }}\left\{
\sum_{j=0}^{\left\vert n\right\vert }M_{j}^{\alpha }\left\vert
K_{M_{j}}\right\vert \right\} .
\end{equation*}
\end{lemma}

\begin{lemma}
\label{Lemma4}Let $0<\alpha \leq 1$ and $\{q_{n}:n\geq 0\}$ be a sequence of
non-increasing numbers, satisfying conditions (\ref{6a}) and (\ref{7a}). If $%
r\geq M_{N}$, then%
\begin{equation*}
\int_{I_{N}}\left\vert F_{r}\left( x-t\right) \right\vert d\mu \left(
t\right) \leq \frac{c_{\alpha }M_{l}^{\alpha }M_{k}}{r^{\alpha }M_{N}},\
x\in I_{l+1}\left( s_{k}e_{k}+s_{l}e_{l}\right) ,
\end{equation*}%
where 
\begin{equation*}
1\leq s_{k}\leq m_{k}-1,\ 1\leq s_{l}\leq m_{l}-1\ \ (k=0,\dots
,N-2,l=k+2,\dots ,N-1)
\end{equation*}%
and 
\begin{equation*}
\int_{I_{N}}\left\vert F_{r}\left( x-t\right) \right\vert d\mu \left(
t\right) \leq \frac{c_{\alpha }M_{k}}{M_{N}},\ x\in I_{N}\left(
s_{k}e_{k}\right) ,
\end{equation*}%
where 
\begin{equation*}
1\leq s_{k}\leq m_{k}-1,(k=0,\dots ,N-1).
\end{equation*}
\end{lemma}

\section{Proofs}

\bigskip

\bigskip \textbf{Proof of Lemma \ref{Lemma2}}. In \cite{go1} Goginava proved
similar equality for the kernel of N\"orlund logarithmic mean $L_{n}$. We
will use his method.

\bigskip Let $sM_{n}<r\leq \left( s+1\right) M_{n},$ where $1\leq s<m_{n}.$
It is easy to show that 
\begin{equation}
\overset{r}{\underset{k=1}{\sum }}q_{r-k}D_{k}=\overset{sM_{n}}{\underset{l=1%
}{\sum }}q_{r-l}D_{l}+\overset{r}{\underset{l=sM_{n}+1}{\sum }}q_{r-l}D_{l}
\label{nor1}
\end{equation}%
\begin{equation*}
:=I+II.
\end{equation*}

By combining (\ref{3dn}) and Abel transformation we get that 
\begin{equation}
I=\overset{sM_{n}-1}{\underset{l=0}{\sum }}q_{r-sM_{n}+l}D_{sM_{n}-l}
\label{nor2}
\end{equation}%
\begin{equation*}
=\overset{sM_{n}-1}{\underset{l=1}{\sum }}%
q_{r-sM_{n}+l}D_{sM_{n}-l}+q_{r-sM_{n}}D_{sM_{n}}
\end{equation*}%
\begin{equation*}
=D_{sM_{n}}\overset{sM_{n}-1}{\underset{l=0}{\sum }}q_{r-sM_{n}+l}
\end{equation*}%
\begin{equation*}
-w_{sM_{n}-1}\overset{sM_{n}-1}{\underset{l=1}{\sum }}q_{r-sM_{n}+l}%
\overline{D_{l}}
\end{equation*}%
\begin{equation}
=\left( Q_{r}-Q_{r-sM_{n}}\right) D_{sM_{n}}  \notag
\end{equation}%
\begin{equation*}
-w_{sM_{n}-1}\overset{sM_{n}-2}{\underset{l=1}{\sum }}\left(
q_{r-sM_{n}+l}-q_{r-sM_{n}+l+1}\right) l\overline{K_{l}}
\end{equation*}%
\begin{equation*}
-w_{sM_{n}-1}q_{r-1}\left( sM_{n}-1\right) \overline{K_{sM_{n}-1}}.
\end{equation*}

Since 
\begin{equation*}
D_{j+sM_{n}}=D_{sM_{n}}+w_{sM_{n}}D_{j},\ \ j=1,2,\dots,sM_{n}-1,
\end{equation*}%
for $II$ we have that 
\begin{equation}
II=\overset{r-sM_{n}}{\underset{l=1}{\sum }}%
q_{r-sM_{n}-l}D_{l+sM_{n}}=Q_{r-sM_{n}}D_{sM_{n}}+w_{sM_{n}}Q_{r-sM_{n}}F_{r-sM_{n}}.
\label{nor3}
\end{equation}

By combining (\ref{nor1})-(\ref{nor3}) we obtain (\ref{8a}) and the proof is
complete.

\textbf{Proof of Lemma \ref{Lemma3}.} Let $sM_{n}<k\leq \left( s+1\right)
M_{n},$ where $1\leq s<m_{n}\ $\ and sequence $\{q_{k}:k\geq 0\}$ be
non-increasing, satisfying condition 
\begin{equation}
\frac{q_{0}n}{Q_{n}}=O(1) ,\text{ \ \ \ \ as \ }n\rightarrow \infty .
\label{100}
\end{equation}

\bigskip By using Abel transformation we get that 
\begin{equation}
Q_{n}=\overset{n}{\underset{j=1}{\sum }}q_{n-j}\cdot 1=\overset{n-1}{%
\underset{j=1}{\sum }}\left( q_{n-j}-q_{n-j-1}\right) j+q_{0}n  \label{2bb}
\end{equation}%
and%
\begin{equation}
F_{n}=\frac{1}{Q_{n}}\left( \overset{n-1}{\underset{j=1}{\sum }}\left(
q_{n-j}-q_{n-j-1}\right) jK_{j}+q_{0}nK_{n}\right) .  \label{2c}
\end{equation}

Since 
\begin{equation}
n\left\vert K_{n}\right\vert \leq c\sum_{A=0}^{\left\vert n\right\vert
}M_{A}\left\vert K_{M_{A}}\right\vert ,  \label{5}
\end{equation}%
by combining (\ref{2bb}) and (\ref{2c}) we immediately get that%
\begin{equation*}
\left\vert F_{n}\right\vert \leq \frac{c}{Q_{n}}\left( \overset{n-1}{%
\underset{j=1}{\sum }}\left\vert q_{n-j}-q_{n-j-1}\right\vert +q_{0}\right)
\sum_{A=0}^{\left\vert n\right\vert }M_{A}\left\vert K_{M_{A}}\right\vert 
\end{equation*}%
\begin{equation*}
=\frac{c}{Q_{n}}\left( \overset{n-1}{\underset{j=1}{\sum }}-\left(
q_{n-j}-q_{n-j-1}\right) +q_{0}\right) \sum_{A=0}^{\left\vert n\right\vert
}M_{A}\left\vert K_{M_{A}}\right\vert 
\end{equation*}%
\begin{equation*}
\leq c\frac{2q_{0}-q_{n-1}}{Q_{n}}\sum_{A=0}^{\left\vert n\right\vert
}M_{A}\left\vert K_{M_{A}}\right\vert 
\end{equation*}%
\begin{equation*}
\leq \frac{c}{Q_{n}}\sum_{A=0}^{\left\vert n\right\vert }M_{A}\left\vert
K_{M_{A}}\right\vert \leq \frac{c}{n}\sum_{A=0}^{\left\vert n\right\vert
}M_{A}\left\vert K_{M_{A}}\right\vert .
\end{equation*}

\bigskip Since the case $q_{0}n/Q_{n}=\overset{\_}{O}\left( 1\right) ,$\ as
\ $n\rightarrow \infty ,$ have already been considered, we can exclude it.

Let $0<\alpha <1.$ We may assume that $\{q_{k}:k\geq 0\}$ satisfies
conditions (\ref{6a}) and (\ref{7a}) and in addition, satisfies the
following 
\begin{equation*}
\frac{Q_{n}}{q_{0}n}=o(1) ,\text{ \ \ \ \ as \ }n\rightarrow \infty .
\end{equation*}

It follows that 
\begin{equation}
q_{n}=q_{0}\frac{q_{n}n}{q_{0}n}\leq q_{0}\frac{Q_{n}}{q_{0}n}=o(1), \text{
\ \ \ \ as \ }n\rightarrow \infty  \label{102}
\end{equation}%
By using (\ref{102}) we immediately get that 
\begin{equation}
q_{n}=\overset{\infty }{\underset{l=n}{\sum }}\left( q_{l}-q_{l+1}\right)
\leq \overset{\infty }{\underset{l=n}{\sum }}\frac{1}{l^{2-\alpha }}\leq 
\frac{c}{n^{1-\alpha }}  \label{103}
\end{equation}%
and%
\begin{equation}
Q_{n}=\overset{n-1}{\underset{l=0}{\sum }}q_{l}\leq \overset{n}{\underset{l=1%
}{\sum }}\frac{c}{l^{1-\alpha }}\leq cn^{\alpha }.  \label{104}
\end{equation}

It is easy to show that 
\begin{equation}
Q_{k}\left\vert D_{sM_{n}}\right\vert \leq cM_{n}^{\alpha }\left\vert
D_{sM_{n}}\right\vert  \label{a1}
\end{equation}%
and 
\begin{equation}
\left( sM_{n}-1\right) q_{k-1}\left\vert K_{sM_{n}-1}\right\vert \leq
ck^{\alpha -1}M_{n}\left\vert K_{sM_{n}-1}\right\vert \leq cM_{n}^{\alpha
}\left\vert K_{sM_{n}-1}\right\vert .  \label{a2}
\end{equation}

Let 
\begin{equation*}
n=s_{n_{1}}M_{n_{1}}+s_{n_{2}}M_{n_{2}}+\dots+ s_{n_{r}}M_{n_{r}},\ \
n_{1}>n_{2}>\dots>n_{r},
\end{equation*}
and 
\begin{equation*}
n^{\left( k\right) }=s_{n_{k+1}}M_{n_{k+1}}+\dots+s_{n_{r}}M_{n_{r}}, \
1\leq s_{n_{l}}\leq m_{l}-1,\ l=1,\dots,r.
\end{equation*}

By combining (\ref{a1}), (\ref{a2}) and Lemma \ref{Lemma2} we have that 
\begin{equation*}
\left\vert Q_{n}F_{n}\right\vert
\end{equation*}%
\begin{equation*}
\leq c_{\alpha }\left( M_{n_{1}}^{\alpha }\left\vert
D_{s_{n_{1}}M_{n_{1}}}\right\vert +\overset{s_{n_{1}}M_{n_{1}}-1}{\underset{%
l=1}{\sum }}\left\vert \left( n^{\left( 1\right) }+l\right) ^{\alpha
-2}\right\vert \left\vert lK_{l}\right\vert +M_{n_{1}}^{\alpha }\left\vert
K_{s_{n_{1}}M_{n_{1}}-1}\right\vert +\left\vert Q_{n^{\left( 1\right)
}}F_{n^{\left( 1\right) }}\right\vert \right) .
\end{equation*}

By repeating this process $r$-times we get that 
\begin{equation*}
\left\vert Q_{n}F_{n}\right\vert
\end{equation*}%
\begin{equation*}
\leq c_{\alpha }\overset{r}{\underset{k=1}{\sum }}\left( M_{n_{k}}^{\alpha
}\left\vert D_{s_{n_{k}}M_{n_{k}}}\right\vert +\overset{s_{n_{k}}M_{n_{k}}-1}%
{\underset{l=1}{\sum }}\left( n^{\left( k\right) }+l\right) ^{\alpha
-2}\left\vert lK_{l}\right\vert +M_{n_{k}}^{\alpha }\left\vert
K_{s_{n_{k}}M_{n_{k}}-1}\right\vert \right)
\end{equation*}%
\begin{equation*}
:=I+II+III.
\end{equation*}

By combining (\ref{1dn}), (\ref{2dn}) and (\ref{5a}) we obtain that 
\begin{equation*}
I\leq c_{\alpha}\overset{\left\vert n\right\vert }{\underset{k=1}{\sum }}%
M_{k}^{\alpha }\left\vert D_{M_{k}}\right\vert \leq c_{\alpha}\overset{%
\left\vert n\right\vert }{\underset{k=1}{\sum }}M_{k}^{\alpha }\left\vert
K_{M_{k}}\right\vert
\end{equation*}%
and 
\begin{equation*}
III\leq c_{\alpha}\overset{r}{\underset{k=1}{\sum }}M_{n_{k}}^{\alpha
-1}\left\vert
M_{n_{k}}K_{s_{n_{k}}M_{n_{k}}}-D_{s_{n_{k}}M_{n_{k}}}\right\vert
\end{equation*}%
\begin{equation*}
\leq c_{\alpha}\overset{r}{\underset{k=1}{\sum }}M_{k}^{\alpha }\left\vert
K_{M_{k}}\right\vert.
\end{equation*}

Moreover, 
\begin{equation*}
II=c_{\alpha }\overset{r}{\underset{k=1}{\sum }}\overset{n_{k}}{\underset{A=1%
}{\sum }}\overset{s_{A}M_{A}-1}{\underset{l=s_{A-1}M_{A-1}}{\sum }}\left(
n^{\left( k\right) }+l\right) ^{\alpha -2}\left\vert lK_{l}\right\vert
\end{equation*}%
\begin{equation*}
=c_{\alpha }\overset{r}{\underset{k=1}{\sum }}\overset{n_{k+1}}{\underset{A=1%
}{\sum }}\overset{s_{A}M_{A}-1}{\underset{l=s_{A-1}M_{A-1}}{\sum }}\left(
n^{\left( k\right) }+l\right) ^{\alpha -2}\left\vert lK_{l}\right\vert
\end{equation*}%
\begin{equation*}
+c_{\alpha }\overset{r}{\underset{k=1}{\sum }}\overset{n_{k}}{\underset{%
A=n_{k+1}+1}{\sum }}\overset{s_{A}M_{A}-1}{\underset{l=s_{A-1}M_{A-1}}{\sum }%
}\left( n^{\left( k\right) }+l\right) ^{\alpha -2}\left\vert
lK_{l}\right\vert
\end{equation*}%
\begin{equation*}
\leq c_{\alpha }\overset{r}{\underset{k=1}{\sum }}M_{n_{k+1}}^{\alpha -2}%
\overset{n_{k+1}}{\underset{A=1}{\sum }}\overset{s_{A}M_{A}-1}{\underset{%
l=s_{A-1}M_{A-1}}{\sum }}\left\vert lK_{l}\right\vert
\end{equation*}%
\begin{equation*}
+c_{\alpha }\overset{r}{\underset{k=1}{\sum }}\overset{n_{k}}{\underset{%
A=n_{k+1}+1}{\sum }}M_{A}^{\alpha -2}\overset{s_{A}M_{A}-1}{\underset{%
l=s_{A-1}M_{A-1}}{\sum }}\left\vert lK_{l}\right\vert
\end{equation*}%
\begin{equation*}
:=II_{1}+II_{2}.
\end{equation*}

By combining (\ref{5a}) and (\ref{5}) for $II_{1}$ we get that 
\begin{equation*}
II_{1}\leq c_{\alpha }\overset{r}{\underset{k=1}{\sum }}M_{n_{k+1}}^{\alpha
-2}\overset{n_{k+1}}{\underset{A=1}{\sum }}\overset{s_{A}M_{A}-1}{\underset{%
l=s_{A-1}M_{A-1}}{\sum }}\overset{A}{\underset{j=0}{\sum }}M_{j}\left\vert
K_{M_{j}}\right\vert
\end{equation*}%
\begin{equation*}
\leq c_{\alpha }\overset{n_{1}}{\underset{k=1}{\sum }}M_{k}^{\alpha -2}%
\overset{k}{\underset{A=1}{\sum }}M_{A}\overset{A}{\underset{j=0}{\sum }}%
M_{j}\left\vert K_{M_{j}}\right\vert
\end{equation*}%
\begin{equation*}
\leq c_{\alpha }\overset{n_{1}}{\underset{k=0}{\sum }}M_{k}^{\alpha -1}%
\overset{k}{\underset{j=0}{\sum }}M_{j}|K_{M_{j}}|
\end{equation*}%
\begin{equation*}
=c_{\alpha }\overset{n_{1}}{\underset{j=0}{\sum }}M_{j}|K_{M_{j}}|\overset{%
n_{1}}{\underset{k=j}{\sum }}M_{k}^{\alpha -1}
\end{equation*}%
\begin{equation*}
\leq c_{\alpha }\overset{n_{1}}{\underset{j=0}{\sum }}M_{j}^{\alpha
}\left\vert K_{M_{j}}\right\vert .
\end{equation*}

\bigskip By using (\ref{5}) for $II_{2}$ we have similarly that 
\begin{equation*}
II_{2}\leq c_{\alpha }\overset{r}{\underset{k=1}{\sum }}\overset{n_{k}}{%
\underset{A=n_{k+1}+1}{\sum }}M_{A}^{\alpha -1}\overset{A}{\underset{j=0}{%
\sum }}M_{j}\left\vert K_{M_{j}}\right\vert
\end{equation*}%
\begin{equation*}
\leq c_{\alpha }\overset{n_{1}}{\underset{A=1}{\sum }}M_{A}^{\alpha -1}%
\overset{A}{\underset{j=0}{\sum }}M_{j}\left\vert K_{M_{j}}\right\vert \leq
c_{\alpha }\overset{n_{1}}{\underset{j=0}{\sum }}M_{j}^{\alpha }\left\vert
K_{M_{j}}\right\vert .
\end{equation*}

The proof is complete by combining the estimates above.

\textbf{Proof of Lemma \ref{Lemma4}.} Let $x\in I_{l+1}\left(
s_{k}e_{k}+s_{l}e_{l}\right),\ 1\leq s_{k}\leq m_{k}-1,\ 1\leq s_{l}\leq
m_{l}-1.$ Then, by applying (\ref{5a}), we have that 
\begin{equation*}
K_{M_{n}}\left( x\right) =0,\text{ when }n>l>k.
\end{equation*}%
Suppose that $k<n\leq l.$ Moreover, by using (\ref{5a}) we get that%
\begin{equation*}
\left\vert K_{M_{n}}\left( x\right) \right\vert \leq cM_{k}.
\end{equation*}%
Let $n\leq k<l.$ Then 
\begin{equation*}
\left\vert K_{M_{n}}\left( x\right) \right\vert =\left( M_{n}+1\right)
/2\leq cM_{k}.
\end{equation*}%
If we now apply Lemma \ref{Lemma3} we can conclude that 
\begin{equation}
Q_{r}\left\vert F_{r}\left( x\right) \right\vert \leq c_{\alpha }\overset{l}{%
\underset{A=0}{\sum }}M_{A}^{\alpha }\left\vert K_{M_{A}}\left( x\right)
\right\vert  \label{7}
\end{equation}%
\begin{equation*}
\leq c_{\alpha }\overset{l}{\underset{A=0}{\sum }}M_{A}^{\alpha }M_{k}\leq
c_{\alpha }M_{l}^{\alpha }M_{k}.
\end{equation*}

Let $x\in I_{l+1}\left( s_{k}e_{k}+s_{l}e_{l}\right),$ for some $0\leq
k<l\leq N-1.$ Since $x-t\in I_{l+1}\left( s_{k}e_{k}+s_{l}e_{l}\right),$ for 
$t\in I_{N}$ and $r\geq M_{N}$ from (\ref{7}) we obtain that 
\begin{equation}
\int_{I_{N}}\left\vert F_{r}\left( x-t\right) \right\vert d\mu \left(
t\right) \leq \frac{c_{\alpha }M_{l}^{\alpha }M_{k}}{r^{\alpha }M_{N}}.
\label{8}
\end{equation}

Let $x\in I_{N}\left( s_{k}e_{k}\right),\ k=0,\dots,N-1.$ Then, by applying
Lemma \ref{Lemma3} and (\ref{5a}) we have that 
\begin{equation}
\int_{I_{N}}Q_{r}\left\vert F_{r}\left( x-t\right) \right\vert d\mu \left(
t\right) \leq \underset{A=0}{\overset{\left\vert r\right\vert }{\sum }}%
M_{A}^{\alpha }\int_{I_{N}}\left\vert K_{M_{A}}\left( x-t\right) \right\vert
d\mu \left( t\right).  \label{9}
\end{equation}

Let $x\in I_{N}\left( s_{k}e_{k}\right),\ k=0,\dots,N-1,\ t\in I_{N}$ and $%
x_{q}\neq t_{q},$ where $N\leq q\leq \left\vert r\right\vert -1.$ By
combining (\ref{5a}) and (\ref{9}) we get that 
\begin{equation*}
\int_{I_{N}}Q_{r}\left\vert F_{r}\left( x-t\right) \right\vert d\mu \left(
t\right)
\end{equation*}%
\begin{equation*}
\leq c_{\alpha }\underset{A=0}{\overset{q-1}{\sum }}M_{A}^{\alpha
}\int_{I_{N}}M_{k}d\mu \left( t\right) \leq \frac{c_{\alpha
}M_{k}M_{q}^{\alpha }}{M_{N}}.
\end{equation*}

Hence, 
\begin{equation}
\int_{I_{N}}\left\vert F_{r}\left( x-t\right) \right\vert d\mu \left(
t\right) \leq \frac{c_{\alpha }M_{k}M_{q}^{\alpha }}{r^{\alpha }M_{N}}\leq 
\frac{c_{\alpha }M_{k}}{M_{N}}.  \label{10}
\end{equation}

Let $x\in I_{N}\left( s_{k}e_{k}\right),\ k=0,\dots,N-1,\ t\in I_{N}$ and $%
x_{N}=t_{N},\dots,x_{\left\vert r\right\vert -1}=t_{\left\vert r\right\vert
-1}.$ By applying again (\ref{5a}) and (\ref{9}) we have that 
\begin{equation}
\int_{I_{N}}\left\vert F_{r}\left( x-t\right) \right\vert d\mu \left(
t\right) \leq \frac{c_{\alpha }}{r^{\alpha }}\overset{\left\vert
r\right\vert -1}{\underset{A=0}{\sum }}M_{A}^{\alpha }\int_{I_{N}}M_{k}d\mu
\left( t\right) \leq \frac{c_{\alpha }M_{k}}{M_{N}}.  \label{11}
\end{equation}

By combining (\ref{8}), (\ref{10}) and (\ref{11}) we complete the proof of \
Lemma \ref{Lemma4}.

\textbf{Proof of Theorem \ref{Theorem1}.} According to Lemma \ref{W} the
proof of the first part of Theorem \ref{Theorem1} will be complete, if we
show that 
\begin{equation*}
\int_{\overline{I_{N}}}\left\vert \overset{\sim }{t}^{\ast }a(x)\right\vert
^{1/\left( 1+\alpha \right) }d\mu \left( x\right) <\infty,
\end{equation*}%
for every $1/\left( 1+\alpha \right) $-atom $a.$ We may assume that $a$ is
an arbitrary $1/\left( 1+\alpha \right) $-atom with support $I,\ \mu \left(
I\right) =M_{N}^{-1}$ and $I=I_{N}.$ It is easy to see that $t_{n}\left(
a\right) =0,$ when $n\leq M_{N}.$ Therefore, we can suppose that $n>M_{N}.$

Let $x\in I_{N}.$ Since $t_{n}$ is bounded from $L_{\infty }$ to $L_{\infty
} $ (the boundedness follows from Lemma 3) and $\left\Vert a\right\Vert
_{\infty }\leq M_{N}^{1+\alpha }$ we obtain that 
\begin{equation*}
\left\vert t_{n}a\left( x\right) \right\vert \leq \int_{I_{N}}\left\vert
a\left( t\right) \right\vert \left\vert F_{n}\left( x-t\right) \right\vert
d\mu \left( t\right)
\end{equation*}%
\begin{equation*}
\leq \left\Vert a\right\Vert _{\infty }\int_{I_{N}}\left\vert F_{n}\left(
x-t\right) \right\vert d\mu \left( t\right)
\end{equation*}%
\begin{equation*}
\leq c_{\alpha }M_{N}^{1+\alpha }\int_{I_{N}}\left\vert F_{n}\left(
x-t\right) \right\vert d\mu \left( t\right) .
\end{equation*}

Let $x\in I_{l+1}\left( s_{k}e_{k}+s_{l}e_{l}\right),\,0\leq k<l<N.$ From
Lemma \ref{Lemma4} we get that 
\begin{equation}
\left\vert t_{n}a\left( x\right) \right\vert \leq \frac{c_{\alpha
}M_{l}^{\alpha }M_{k}M_{N}^{\alpha }}{n^{\alpha }}.  \label{12}
\end{equation}

Let $x\in I_{N}\left( s_{k}e_{k}\right),\,0\leq k<N.$ From Lemma \ref{Lemma4}
we have that 
\begin{equation}
\left\vert t_{n}a\left( x\right) \right\vert \leq c_{\alpha
}M_{k}M_{N}^{\alpha }.  \label{12a}
\end{equation}

By combining (\ref{2}) and (\ref{12})-(\ref{12a}) we obtain that 
\begin{equation*}
\int_{\overline{I_{N}}}\left\vert \overset{\sim }{t^{\ast }}a(x)\right\vert
^{1/\left( 1+\alpha \right) }d\mu \left( x\right)
\end{equation*}%
\begin{equation*}
=\overset{N-2}{\underset{k=0}{\sum }}\overset{m_{k}-1}{\underset{s_{k}=1}{%
\sum }}\overset{N-1}{\underset{l=k+1}{\sum }}\overset{m_{l}-1}{\underset{%
s_{l}=1}{\sum }}\int_{I_{l+1}\left( s_{k}e_{k}+s_{l}e_{l}\right) }\underset{%
n>M_{N}}{\sup }\left\vert \frac{t_{n}a\left( x\right) }{\log ^{1+\alpha
}(n+1)}\right\vert ^{1/\left( 1+\alpha \right) }d\mu \left( x\right)
\end{equation*}%
\begin{equation*}
+\overset{N-1}{\underset{k=0}{\sum }}\overset{m_{k}-1}{\underset{s_{k}=1}{%
\sum }}\int_{I_{N}\left( s_{k}e_{k}\right) }\underset{n>M_{N}}{\sup }%
\left\vert \frac{t_{n}a\left( x\right) }{\log ^{1+\alpha }(n+1)}\right\vert
^{1/\left( 1+\alpha \right) }d\mu \left( x\right)
\end{equation*}%
\begin{equation*}
\leq \frac{c_{\alpha }}{N}\overset{N-2}{\underset{k=0}{\sum }}\overset{%
m_{k}-1}{\underset{s_{k}=1}{\sum }}\overset{N-1}{\underset{l=k+1}{\sum }}%
\overset{m_{l}-1}{\underset{s_{l}=1}{\sum }}\int_{I_{l+1}\left(
s_{k}e_{k}+s_{l}e_{l}\right) }\underset{n>M_{N}}{\sup }\left\vert
t_{n}a\left( x\right) \right\vert ^{1/\left( 1+\alpha \right) }d\mu \left(
x\right)
\end{equation*}%
\begin{equation*}
+\frac{c_{\alpha }}{N}\overset{N-1}{\underset{k=0}{\sum }}\overset{m_{k}-1}{%
\underset{s_{k}=1}{\sum }}\int_{I_{N}\left( s_{k}e_{k}\right) }\underset{%
n>M_{N}}{\sup }\left\vert t_{n}a\left( x\right) \right\vert ^{1/\left(
1+\alpha \right) }d\mu \left( x\right)
\end{equation*}%
\begin{equation*}
\leq \frac{c_{\alpha }}{N}\overset{N-2}{\underset{k=0}{\sum }}\overset{N-1}{%
\underset{l=k+1}{\sum }}\frac{\left( m_{k}-1\right) \left( m_{l}-1\right) }{%
M_{l+1}}\frac{\left( M_{l}^{\alpha }M_{k}\right) ^{1/\left( 1+\alpha \right)
}M_{N}^{\alpha /\left( 1+\alpha \right) }}{n^{\alpha /\left( 1+\alpha
\right) }}
\end{equation*}%
\begin{equation*}
+\frac{c_{\alpha }}{N}\overset{N-1}{\underset{k=0}{\sum }}\frac{\left(
m_{k}-1\right) }{M_{N}}M_{N}^{\alpha /\left( 1+\alpha \right)
}M_{k}^{1/\left( 1+\alpha \right) }
\end{equation*}%
\begin{equation*}
\leq \frac{c_{\alpha }M_{N}^{\alpha /\left( 1+\alpha \right) }}{Nn^{\alpha
/\left( 1+\alpha \right) }}\overset{N-2}{\underset{k=0}{\sum }}\overset{N-1}{%
\underset{l=k+1}{\sum }}\frac{\left( M_{l}^{\alpha }M_{k}\right) ^{1/\left(
1+\alpha \right) }}{M_{l+1}}
\end{equation*}%
\begin{equation*}
+\frac{c_{\alpha }}{N}\overset{N-1}{\underset{k=0}{\sum }}\frac{%
M_{k}^{1/\left( 1+\alpha \right) }}{M_{N}^{1/\left( 1+\alpha \right) }}\leq
c_{\alpha }<\infty .
\end{equation*}%
The proof is complete.

\textbf{Proof of Theorem \ref{Theorem2}. }By Lemma \ref{W} the proof of
Theorem \ref{Theorem2} will be complete, if we show that

\begin{equation*}
\frac{1}{\log n}\overset{n}{\underset{k=1}{\sum }}\frac{\left\Vert
t_{k}a\right\Vert _{1/\left( 1+\alpha \right) }^{1/\left( 1+\alpha \right) }%
}{k}\leq c_{\alpha }<\infty,
\end{equation*}%
for every $1/\left( 1+\alpha \right) $-atom $a.$ Analogously to the proof of
Theorem \ref{Theorem1} we may assume that $a$ be an arbitrary $1/\left(
1+\alpha \right) $-atom with support $I,\ \mu \left( I\right) = M_{N}^{-1}$
and $I=I_{N}$ and $n>M_{N}.$

Let $x\in I_{N}.$ Since $t_{m}$ is bounded from $L_{\infty }$ to $L_{\infty
} $ (the boundedness follows from Lemma 3) and $\left\Vert a\right\Vert
_{\infty }\leq M_{N}^{1+\alpha }$, we obtain that 
\begin{equation*}
\int_{I_{N}}\left\vert t_{n}a\left( x\right) \right\vert ^{1/\left( 1+\alpha
\right) }d\mu \leq \left\Vert a\left( x\right) \right\Vert _{\infty
}^{1/\left( 1+\alpha \right) }M_{N}^{-1}\leq c_{\alpha }<\infty .
\end{equation*}%
Hence 
\begin{equation*}
\frac{1}{\log n}\overset{n}{\underset{k=M_{N}}{\sum }}\frac{%
\int_{I_{N}}\left\vert t_{k}a\left( x\right) \right\vert ^{1/\left( 1+\alpha
\right) }d\mu }{k}
\end{equation*}%
\begin{equation*}
\leq \frac{c_{\alpha }}{\log n}\overset{n}{\underset{k=1}{\sum }}\frac{1}{k}%
\leq c_{\alpha }<\infty .
\end{equation*}

By combining (\ref{2}) and (\ref{12})-(\ref{12a}) we can conclude that 
\begin{equation*}
\frac{1}{\log n}\overset{n}{\underset{k=M_{N}+1}{\sum }}\frac{\int_{%
\overline{I_{N}}}\left\vert t_{k}a\left( x\right) \right\vert ^{1/\left(
1+\alpha \right) }d\mu \left( x\right) }{k}
\end{equation*}%
\begin{equation*}
=\frac{1}{\log n}\overset{n}{\underset{k=M_{N}+1}{\sum }}\overset{N-2}{%
\underset{r=0}{\sum }}\overset{m_{r}-1}{\underset{s_{r}=1}{\sum }}\overset{%
N-1}{\underset{l=r+1}{\sum }}\overset{m_{l}-1}{\underset{s_{l}=1}{\sum }}%
\frac{\int_{I_{l+1}\left( s_{r}e_{r}+s_{l}e_{l}\right) }\left\vert
t_{k}a\left( x\right) \right\vert ^{1/\left( 1+\alpha \right) }d\mu \left(
x\right) }{k}
\end{equation*}%
\begin{equation*}
+\frac{1}{\log n}\overset{n}{\underset{k=M_{N}+1}{\sum }}\overset{N-1}{%
\underset{r=0}{\sum }}\overset{m_{r}-1}{\underset{s_{r}=1}{\sum }}\frac{%
\int_{I_{N}\left( s_{r}e_{r}\right) }\left\vert t_{k}a\left( x\right)
\right\vert ^{1/\left( 1+\alpha \right) }d\mu \left( x\right) }{k}
\end{equation*}%
\begin{equation*}
\leq \frac{1}{\log n}\left( \overset{n}{\underset{k=M_{N}+1}{\sum }}\frac{%
c_{\alpha }M_{N}^{\alpha /\left( 1+\alpha \right) }}{k^{\alpha /\left(
1+\alpha \right) +1}}+\overset{n}{\underset{k=M_{N}+1}{\sum }}\frac{%
c_{\alpha }}{k}\right) <c_{\alpha }<\infty .
\end{equation*}%
The proof is complete.

\textbf{Proof of Remark \ref{Remark1}}. Let us see an example. Let 
\begin{equation*}
q_{n}:=\left\{ 
\begin{array}{ll}
\frac{1}{\sqrt{n}} & \text{ if }n\in \mathbb{N}_{+} \\ 
0 & \text{ if }n=0.%
\end{array}%
\right.
\end{equation*}%
Then sequence is non-increasing and non-negative.

1. Let $0<\alpha <1/2$ arbitrary. It is easy to see, that if $n>2$, then 
\begin{equation*}
Q_{n}>\sum_{k=1}^{n-2}\frac{1}{\sqrt{k}}>\int_{1}^{n-1}\frac{1}{\sqrt{x}}
dx=2\sqrt{n-1}-2
\end{equation*}

Recalling $\alpha <1/2$ we obtain 
\begin{equation*}
0<\frac{n^{\alpha }}{Q_{n}}<\frac{n^{\alpha }}{2\sqrt{n-1}-2}=O(1)
\end{equation*}

On the other hand 
\begin{equation*}
\frac{q_{n}-q_{n+1}}{n^{\alpha -2}}=\frac{1}{\sqrt{1+\frac{1}{n}}\left( 
\sqrt{1+\frac{1}{n}}+1\right) }n^{\frac{1}{2}-\alpha }\neq O(1).
\end{equation*}

2. Analogously we can show that in the case of $\alpha >1/2$ the situation
is the opposite.

\end{document}